\input amstex
\documentstyle{amsppt}
\magnification 1000
\catcode`\@=11
\def\logo\@{}
\catcode`\@=\active  
\NoBlackBoxes
\baselineskip=22pt
\topmatter
\title Comment on "Existence and Regularity  for an Energy Maximization Problem
in Two Dimensions" \endtitle
\author 
\endauthor
\endtopmatter
\document

\define\[{\left[}%
\define\]{\right]}%
\define\({\left(}%
\define\){\right)}%
\baselineskip=20pt

Spyridon Kamvissis

Max Planck Institute for
Mathematics in the Sciences, Leipzig, Germany

and

Department of Applied Mathematics, University of Crete, Greece

\bigskip

ABSTRACT

A revision of the last appendix of the paper 
"Existence and Regularity  for an Energy Maximization Problem
in Two Dimensions" by S.Kamvissis and E.A.Rakhmanov,
that appeared in the Journal of Mathematical Physics,  v.46, n.8, 2005.

\newpage

DROPPING ASSUMPTION (A) IN SECTION 5 OF [9].

\bigskip

In section 5 of [9], we have assumed that the solution of the problem of the
maximization 
of the equilibrium energy  is a continuum, say $F$, which does not intersect the
linear
segment $[0, iA]$ except of course  at $0_+, 0_-$.
We also prove that $F$ does not touch the real line, except of course at
$0$ and possibly $\infty$. This enables us to take variations in section 6,
keeping fixed a finite number of points, and thus arrive at the identity of
Theorem 5,
from which we derive the regularity of $F$ and the fact that $F$ is, after all,
an S-curve. 

In general, it is conceivable that $F$ intersects the linear
segment $[0, iA]$ at points other than  $0_+, 0_-$.
If the set of such points is finite, there is no problem, since we can always
consider
variations keeping fixed a finite number of points, and arrive at the same
result (see the remark after the proof of Theorem 5).

If, on the other hand, this is not the case, we have a different kind of
problem, 
because the function $V$ introduced in section 6 (the complexification of
the field) is not analytic across the segment $[-iA, iA]$.

What is true, however, is that $V$ is analytic in a Riemann surface consisting
of infinitely
many sheets, cut along the line segment $[-iA, iA]$. So, the appropriate, 
underlying space for the 
(doubled up) variational problem should now be a
non-compact Riemann surface, say $\Bbb L$.

Compactness is crucial in the proof of a maximizing continuum. But we can
compactify the
Riemann surface $\Bbb L$ by compactifying the
complex plane.
Let  the map $\Bbb C \to \Bbb L$ be defined by
$$
\aligned
y = log(z-iA) - log(z+iA).
\endaligned
$$
The point $z =iA$ corresponds to infinitely many y-points, i.e.
$y = -\infty + i \theta,~~\theta \in \Bbb R$, which will be identified.
Similarly,
the point $z= -iA$  corresponds to infinitely many points 
$y = +\infty + i \theta,~~\theta \in \Bbb R$, which will also be identified. 
The point $0 \in \Bbb C$  corresponds to the points $ k \pi i$, $k$ odd.

By compactifying the plane we then compactify the Riemann surface $\Bbb L$.
The distance between two points in the Riemann surface $\Bbb L$
is defined to be the corresponding stereographic distance between the
images of these points in the compactified $\Bbb C$.

With these changes, the proof of the existence of the maximizing continuum in
sections 1, 3, 4
goes through virtually unaltered. In section 6, we would have to
consider the complex field $V$ as a function defined in the
Riemann surface $\Bbb L$ and all proofs go through.
The corresponding
result of section 7
will give us an S-curve $C$ in the Riemann surface $\Bbb L$. 
We then have the following facts.

Consider the  image $\Bbb D$ of the closed upper half-plane 
under
$$
\aligned
y = log(z-iA) - log(z+iA).
\endaligned
$$
Consider continua in $\Bbb D$ containing the points $y=\pi i$ and $y=-\pi i$.
Define the Green's potential  and Green's energy of a Borel measure by 
(4), (5), (6) and the equilibrium measure by (7).
Then there exists a continuum
$F$  maximizing the equilibrium energy, for the field given by
(3) with conditions (1). $F$ does not touch $\partial \Bbb D$ except at a finite
number of points. By taking  variations as in section 6, one sees 
that $F$ is an S-curve. In particular,
the support of the equilibrium measure on $F$ is a union of analytic
arcs and 
at any interior point of $supp\mu$
$$
\aligned
{d \over {d n_+}} (\phi + V^{\lambda^F}) =
{d \over {d n_-}}  (\phi + V^{\lambda^F}),
\endaligned
\tag8
$$
where the two derivatives above denote the normal derivatives.

\bigskip

We then have the following.

THEOREM 9. Consider the semiclassical limit ($\hbar \to 0$)
of the solution of (9)-(10) (that is the initial value problem for the
focusing NLS with  parameter $\hbar$)
with bell-shaped initial data. Replace the initial data by the so-called
soliton ensembles data (as introduced in [3])  defined by replacing
the scattering data for $\psi (x,0)= \psi_0(x)$ by their
WKB-approximation.  Assume, for simplicity,
that the spectral density of eigenvalues
satisfies conditions (1).

Then, asymptotically as $\hbar \to 0$, the
solution  $\psi (x,t)$ admits a "finite genus description", in the sense of
Theorem A.1.

PROOF: (i) The proof of the existence of an S-curve $F$ in $\Bbb L$ follows as above.

(ii) We want to deform the original discrete Riemann-Hilbert problem to the set
$\hat F$ consisting of the  projection of  $F$
to the complex plane.
It is clear however that  $\hat F$ may not encircle the spike $[0,iA]$.
It is possible, on the other hand,  to append S-loops (not necessarily with respect to the same
branch of the external field)
and end up with a sum of S-loops, such that the amended $\hat F$  $does$
encircle the spike $[0,iA]$, meaning that  $[0, iA]$ is a subset
of the closure of the union   of the interiors of the loops
of which $ \hat F$ consists. A little thought shows that this is all we need.
(Indeed, within each of the loops we use the same pole-removing transformation as 
in [3]. Eventually of course we will use different interpolations, according to the sheet
of each piece of $F$.)

To see that we can always append the needed S-loop,  suppose
there is an open interval, say $(i\alpha, i\alpha_1)$,
which lies in the exterior of $\hat F$, while
$i\alpha, i\alpha_1 \in \hat F$. Let us assume  for example that $\hat F$ crosses
$[0,iA]$ along bands at $i\alpha, i\alpha_1$; call these bands $S, S_1$.
Let $\beta^-, \beta^+$ be points (considered in $\Bbb C$)
lying on $S$ to the left and right
of $i\alpha$ respectively, and at a small distance
from $i\alpha$. Similarly,
let $\beta_1^-, \beta_1^+$ be points lying on $S_{1}$
to the left and right
of $i\alpha_{1}$ respectively, and at a small distance
from $i\alpha_{1}$.
We will show that there exists a "gap" region including the preimages of
$\beta^-, \beta_{1}^-$ lying in the $N$th sheet for $-N$ large enough,
and similarly there exists a "gap" region including the preimages of
$\beta^+, \beta_{1}^+$ lying in the $M$th sheet for $M$ large enough,
both being  regions for which the gap inequalities hold a priori,
irrespectively of the actual S-curve,
depending only on the external field!

Indeed, note  that the
quantity  $Re (\tilde \phi^{\sigma} (z) )$ (which defines the variational
inequalities) is a priori bounded above  by
$-\phi(z)$. For this, see (8.8) in Chapter 8 of [3]; there is actually a sign error:
the right formula is
$$
\aligned
Re (\tilde \phi^{\sigma} (z) )= -\phi(z)
+ \int G(z,\eta) \rho^{\sigma} (\eta)d\eta.
\endaligned
$$
Next note that the difference of the values of the
function $Re (\tilde \phi^{\sigma} (z) )$ in consecutive sheets is
$\delta Re (\tilde \phi^{\sigma}) = \pm 2 \pi Rez$, and hence the difference of
the values at points on consecutive sheets whose image under the
projection to the complex plane is $i \eta + \epsilon$, where $\eta $ is real
and $\epsilon $ is a small (negative or positive) real,
is $\delta (Re \tilde \phi^{\sigma} )= \pm 2 \pi \epsilon$.
This means that on the left (respectively right)
side of the imaginary semiaxis, the inequality $Re ( \tilde \phi^{\sigma} (z) )<0$
will be eventually (depending
on the sheet) be valid at any given small distance to it.

We now connect  the preimages of $\beta^-$
and  $\beta_{1}^-$ (under the projection of $\Bbb L $ to
$\Bbb C$) lying in the $N$th sheet to 
$\beta^-$ and $\beta_1^-$ respectively, using the
results of [9]. Similarly  we join  
the preimages of $\beta^+$
and  $\beta_{1}^+$ lying in the $M$th sheet
to $\beta^+$ and $\beta_1^+$ respectively.

Then, we  join the  the preimages of $\beta^-$
and  $\beta_{1}^-$ (under the projection of $\Bbb L $ to
$\Bbb C$) lying in the $N$th sheet and the preimages of $\beta^+$
and  $\beta_{1}^+$ lying in the $M$th sheet, along the
according gap regions.

It is easy to see that we  end up with an S-loop whose projection is
covering the "lacuna" $(i\alpha, i\alpha_{1})$.

The original discrete Riemann-Hilbert problem  can be  trivially deformed to a
discrete Riemann-Hilbert on the resulting  (projection of the)
union of S-loops. All this is possible even in the
case where $\hat F$ self-intersects. 

(iii) We  deform the discrete Riemann-Hilbert problem to the continuous one with the right
band/gap structure
(on $\hat F$; according to the  equilibrium measure on $F$), which is
then explicitly solvable via theta functions exactly as in [3].
Both the discrete-to-continuous approximation and the opening of the lenses needed
for this deformation are justified  as in [3] 
and therefore the technical details will not be repeated here.
It is important to notice that our construction has ensured the analytic continuation of
the jump matrix along $\hat F$ (oriented according to $F$).
The g-function is defined by the same Thouless-type formula with respect to the
equilibrium measure (cf. section 2(iii)).
It satisfies the same conditions as in [3]
(measure reality and variational inequality) on bands  and  gaps.
The equilibrium measure lives in $\Bbb L$ but the  Riemann-Hilbert problem lives in $\Bbb C$.

\newpage

REFERENCES (including those appeared in [9])
 
\bigskip

[1]  A. A. Gonchar and E. A. Rakhmanov, Equilibrium Distributions and Degree
of Rational Approximation of Analytic Functions,
Math. USSR Sbornik, v. 62, pp.305--348, 1989.

[2] E. B. Saff, V. Totik, Logarithmic Potentials
with External Fields, Springer Verlag, 1997.

[3] S. Kamvissis, K. T.-R. McLaughlin, P. D. Miller, Semiclassical Soliton
Ensembles
for the Focusing Nonlinear Schr\"odinger Equation,
Annals of Mathematics Studies, v.154, Princeton University Press, 2003.

[4] P. Deift, X.Zhou, A Steepest Descent Method for
Oscillatory Riemann-Hilbert Problems, Annals of Mathematics, v.137, n.2, 1993,
pp.295-368.

[5] E. A. Perevozhnikova and E. A. Rakhmanov, Variations of the Equilibrium
Energy and S-property of Compacta of Minimal Capacity, preprint, 1994.

[6] J. Dieudonn\'e, Foundations of Modern Analysis, Academic Press, 1969.

[7] N. S. Landkof, Foundations of Modern Potential Theory, Springer Verlag,
1972.

[8] G. M. Goluzin, Geometric Theory of Functions of a Complex Variable,
Translations of Mathematical Monographs, v.26, AMS 1969.

[9] S.Kamvissis, E.A.Rakhmanov, Journal of Mathematical Physics,  v.46, n.8, 2005.

\end{document}